\documentclass[11pt]{amsart}

\usepackage{amsfonts}
\usepackage{graphicx}
\usepackage{amsmath}
\usepackage{amssymb}

 \advance\hoffset by -0.0 truecm

\newtheorem{theorem}{Theorem}
\newtheorem{corollary}{Corollary}
\newtheorem{lemma}{Lemma}
\newtheorem{proposition}{Proposition}
\theoremstyle{definition}
\newtheorem{definition}{Definition}
\theoremstyle{remark}
\newtheorem{remark}{Remark}
\numberwithin{equation}{section}

\begin{document}

\title[Nonholonomic Lorentzian geometry]
{Nonholonomic Lorentzian geometry on some $\mathbb H$-type groups}

\author{Anna Korolko,\ Irina Markina}

\address{Department of Mathematics,
University of Bergen, Johannes Brunsgate 12, Bergen 5008, Norway}

\email{anna.korolko@uib.no}

\thanks{}

\address{Department of Mathematics,
University of Bergen, Johannes Brunsgate 12, Bergen 5008, Norway}

\email{irina.markina@uib.no}

\thanks{The authors are supported by a research
grant of the Norwegian Research Council \#177355/V30, and by the European Science Foundation Research Networking Programme HCAA}





\begin{abstract}
We consider examples of the $\mathbb H$-type groups with the natural horizontal distribution generated by the commutation relations of the group. In the contrast with the previous studies we furnish the horizontal distribution with the Lorentzian metric, which is nondegenerate metric of index 1 instead of a positive definite quadratic form. The causal character is defined. We study the reachable set by timelike future directed curves. The parametric equations of geodesics are obtained.

\end{abstract}

\maketitle

\section{Introduction}

Sub-Riemannian manifolds and the geometry introduced by bracket generating distributions of smoothly varying $k$-plains is widely studied, interesting subject, which has applications in control theory, quantum physics, C-R geometry, and other areas. The main difference of the sub-Riemannian manifold from a Riemannian one is the presence of a smooth subbundle of the tangent bundle, generating the entire tangent bundle by means of the commutators of vector fields. The subbundle, which  is often called horizontal, is equipped with a positively definite metric that leads to the triple: manifold, horizontal subbundle, and Riemannian metric on the horizontal subbundle, which is called a {\it sub-Riemannian manifold}. The foundation of the sub-Riemannian geometry can be found in~\cite{LS,Mon2,Mon,Str}. The following question can be asked. What kind of geometrical features will have the mentioned triplet if we change the positively definite metric to an indefinite nondegenerate metric. It is natural to start with the Lorentzian metric of index 1. In this case the triplet: manifold, horizontal subbundle, and Lorentzian metric on the horizontal subbundle can be called {\it sub-Lorentzian manifold}. It was mentioned in~\cite{Str} that it would be interesting to consider the sub-Lorentzian geometry, but, as it is known to the authors, there are only few works devoted to this subject~\cite{ChMV,Groch1,Groch2}. In the present paper we study two examples of $\mathbb H$-type groups furnished with the Lorentzian metric. The notion of $\mathbb H$-type group was introduced in~\cite{Kap1}, see also a nice classification in~\cite{CDKR}. The simplest example of a nonabelian $\mathbb H$-type group equipped with the Lorentzian metric is the Heisenberg group with the Lorentzian metric and it has been considered in~\cite{Groch1,Groch2}. Our example is an extension of the Heisenberg group having $3$-dimensional center related with the quaternion division algebra. The so-called quaternion $\mathbb H$-type group with the positive definite metric was studied in~\cite{ChM1,ChM2,ChM3}. We study the Heisenberg group and the quaternion $\mathbb H$-type group endowed with the Lorentzian metric defining a causal character of the manifold under consideration. We give a description of the reachable set by timelike future directed curves. We find the parametric equations for geodesics on the Heisenberg and quaternion $\mathbb H$-groups with the Lorentzian metric. Unlikely to the sub-Riemannian Heisenberg group, we get the uniqueness of geodesics starting from the origin. Notice that by a geodesic we mean a projection of a solution of the corresponding Hamiltonian system onto the underlying manifold.

\section{Basic definitions and notations}

$\mathbb{H}$-type homogeneous groups are simply connected 2-step Lie groups $\mathbb{G}$ whose algebras $\mathcal{G}$ are graded and carry an inner product such that
\begin{itemize}
 \item [$(i)$]{ $\mathcal{G}$ is the orthogonal direct sum of the generating subspace $V_1$ and the center $V_2\colon \,\mathcal{G}=
V_1\oplus V_2$, $\,V_2=[V_1,V_1]$, $\,[V_1,V_2]=0$,}
\item[$(ii)$]{ the homomorphisms $J_Z\colon V_1\to V_1$, $Z\in V_2$, defined by
\begin{gather}
    \langle J_ZX,X'\rangle=\langle Z,[X,X']\rangle,\quad X,\,X'\in V_1,\notag
\end{gather}
satisfy the equation $J_Z^2=-|Z|^2I$, $\,Z\in V_2$.}
\end{itemize} Here $\langle\cdot,\cdot\rangle$ is an inner product on $\mathcal{G}$, $[\cdot,\cdot]$ is a commutator and $I$ is the identity.

We will consider the quaternion $\mathbb{H}$-type group $\textbf{H}$ with $V_1$ associated with the space of quaternions and $V_2$ as a three-dimensional center. Remind that quaternion numbers are a noncommutative extension of complex numbers. A quaternion can be represented in a matrix way by $q=a\mathcal{U}+b\mathcal{I}_1+c\mathcal{I}_2+d\mathcal{I}_3$, where

\begin{gather} \mathcal{U}=
\begin{bmatrix}
1& 0& 0& 0\\
0& 1& 0& 0\\
0& 0& 1& 0\\
0& 0& 0& 1
\end{bmatrix},
\quad
\mathcal{I}_1=
\begin{bmatrix}
0& 1& 0& 0\\
-1& 0& 0& 0\\
0& 0& 0& 1\\
0& 0& -1& 0
\end{bmatrix},\notag\\
\mathcal{I}_2=
\begin{bmatrix}
0& 0& 0& -1\\
0& 0& -1& 0\\
0& 1& 0& 0\\
1& 0& 0& 0
\end{bmatrix},
\quad
\mathcal{I}_3=
\begin{bmatrix}
0& 0& -1& 0\\
0& 0& 0& 1\\
1& 0& 0& 0\\
0& -1& 0& 0
\end{bmatrix}\notag
\end{gather}
is the basis of quaternion numbers given by real matrices. As a vector space $Q$ over the real numbers, the quaternions have dimension 4, whereas the complex numbers have dimension 2.

$\textbf{H}$ is a set $\mathbb R^4\times \mathbb{R}^3$ with the multiplication law defined in the following way
\begin{gather}
 L_q(q')=L_{(x,z)}(x',z')=(x,z_{1},z_{2},z_{3})\circ (x',z'_{1},z'_{2},z'_{3})\notag\\
=(x+x',z_{1}+z'_{1}+\frac{1}{2}({\mathcal{I}_1}x,x'),\,z_{2}+z'_{2}+\frac{1}{2}({\mathcal{I}_2}x,x'),\,z_{3}+z'_{3}+\frac{1}{2}({\mathcal{I}_3}x,x')),\notag
\end{gather} for $q=(x,z)$ and $q'=(x',z')$, where $x,x'\in \mathbb{R}^4$ and $({\mathcal{I}_1}x,x')$, $({\mathcal{I}_2}x,x')$, $({\mathcal{I}_3}x,x')$ are the usual inner products of the vectors ${\mathcal{I}_1}x,\,{\mathcal{I}_2}x,\,{\mathcal{I}_3}x$ belonging to $\mathbb{R}^4$ by $x'\in \mathbb{R}^4$. The unit element $e$ of the group $\mathbf H$ is $e=(0,0)$ and the inverse element to $(x,z)$ is $(-x,-z)$.
The multiplication ``$\circ$`` defines the left translation $L_q$ of $q'$ by the element $q=(x,z)\in \textbf{H}$ on the group $\textbf{H}$. The topological dimension of the group $\textbf{H}$ is 7, the Hausdorff dimension is 10.
The matrices $\mathcal{I}_1$, $\mathcal{I}_2$ and $\mathcal{I}_3$ define the homomorphisms $J_Z$.

$\mathcal{G}$ is a Lie algebra of the group $\textbf{H}$ which can be associated with the set of all left invariant vector fields of the tangent space $T\textbf{H}_e$ at the unity $e$. $T_h\textbf{H}\subset T\textbf{H}$ is a subbundle of the tangent bundle spanned by the left invariant vector fields $X_{\alpha}(x,z)$ with $X_{\alpha}(0,0)=\frac{\partial }{\partial x_{\alpha}}$, $\alpha=1,\ldots,4$. In coordinates of the standard Euclidean basis $\frac{\partial }{\partial x_{\alpha}}$, $\alpha=1,\ldots,4$, $\frac{\partial}{\partial z_{\beta}}$,
$\beta=1,2,3$, these vector fields are expressed as
\begin{gather}
   X_1=\frac{\partial}{\partial x_1}+\frac{\displaystyle1}{\displaystyle2} \left(+x_2
   \frac{\partial}{\partial z_{1}}-x_4\frac{\partial}{\partial z_{2}}
   - x_3\frac{\partial}{\partial z_{3}}\right),
   \notag\\
   X_2=\frac{\partial}{\partial x_2}+\frac{\displaystyle1}{\displaystyle2} \left(-x_1
   \frac{\partial}{\partial z_{1}}-x_3\frac{\partial}{\partial z_{2}}
   + x_4\frac{\partial}{\partial z_{3}}\right),
   \notag\\
   X_3=\frac{\partial}{\partial x_3}+\frac{\displaystyle1}{\displaystyle2} \left(+x_4
   \frac{\partial}{\partial z_{1}}+x_2\frac{\partial}{\partial z_{2}}
   + x_1\frac{\partial}{\partial z_{3}}\right),
   \notag\\
   X_4=\frac{\partial}{\partial x_4}+\frac{\displaystyle1}{\displaystyle2} \left(-x_3
   \frac{\partial}{\partial z_{1}}+x_1\frac{\partial}{\partial z_{2}}
   - x_2\frac{\partial}{\partial z_{3}}\right).
   \notag
\end{gather} We call $T_h\mathbf{H}$ the {\it horizontal bundle} and any vector field $V\in T_h\mathbf{H}$ is called the horizontal vector field.
The left invariant vector fields $Z_{\beta}=\frac{\partial}{\partial z_{\beta}}$, $\beta=1,2,3$ form a basis of the complement to~$T_h\textbf{H}$.
The commutation relations are given by
\begin{gather}
    [X_1,X_2]=-Z_1,\quad [X_1,X_3]=Z_3,\quad [X_1,X_4]=Z_2,\notag\\
    [X_2,X_3]=Z_2,\quad [X_2,X_4]=-Z_3,\quad [X_3,X_4]=-Z_1.\notag
\end{gather} We see that $X_{\alpha}$, $\alpha=1,\ldots,4$, and their commutators $Z_{\beta}$, $\beta=1,2,3$, span the entire tangent bundle $T\mathbf H$. Therefore, the horizontal bundle $T_h\mathbf H$ is {\it bracket generating} of step $2$, see the definition of bracket generating subbundle in~\cite{Chow, Hormander}.

A basis of one-forms dual to $X_1,\ldots,X_4,Z_1,Z_2,Z_3$ is given by $dx_1,\ldots,dx_4,\upsilon_1,\upsilon_2,\upsilon_3$ with
\begin{gather}
    \upsilon_1=dz_1-\frac{1}{2}(+x_2 dx_1-x_1dx_2+x_4dx_3-x_3dx_4),\notag\\
    \upsilon_2=dz_2-\frac{1}{2}(-x_4 dx_1-x_3dx_2+x_2dx_3+x_1dx_4),\notag\\
    \upsilon_3=dz_3-\frac{1}{2}(-x_3 dx_1+x_4dx_2+x_1dx_3-x_2dx_4).\notag
\end{gather}

We use the normal coordinates $(x,z)=(x_1,\ldots,x_4,z_{1},z_{2},z_{3})$ for the elements
\begin{gather}
 \exp\left(\sum\limits_{\alpha=1}^4x_{\alpha}X_{\alpha}
+\sum\limits_{\beta=1}^3z_{\beta}Z_{\beta}\right)\in \textbf{H}.\notag
\end{gather}

An absolutely continuous curve $c(s)\colon [0,1]\to \textbf{H}$ is called the {\it horizontal curve} if its tangent vector $\dot{c}(s)$ belongs to $V_1$ at each point $c(s)$. Any two points in $\mathbf H$ can be connected by piecewise smooth horizontal curve because of bracket generating property of $T_h\mathbf H$.
Unless otherwise stated all vectors and curves are supposed to be horizontal.

The geometry of $\mathbb H$-type groups was studied, for instance, in~\cite{BGGr2,CChGr3,ChM4,ChM1,ChM3,CDKR,Kap2,Kor1}. In all these cases the horizontal bundle $T_h\mathbf H$ was endowed with a positively definite metric. We will consider the same type of the horizontal bundle but equipped with the Lorentzian metric.
We call such a manifold {\it sub-Lorentzian}, analogously to the sub-Riemannian case, and we formulate it in the following definition.

\begin{definition}
Let $M$ be a smooth $n$-dimensional manifold, let $T_hM$ be a smooth $k$-dimensional, $k<n$, bracket generating subbundle on $TM$, and let$Q_{ T_hM}(\cdot,\cdot)$ be a smooth Lorentzian metric on $T_hM$. Then the triple $(M,T_hM,Q_{T_hM}(\cdot,\cdot))$ is called the sub-Lorentzian manifold.
\end{definition}
We will call $Q_{T_hM}(\cdot,\cdot)$ the {\it sub-Lorentzian metric} and skipping the subscript, we write $Q(\cdot,\cdot)$.

We define a smooth Lorentzian metric on $T_h\mathbf H$ by
\begin{equation}\label{timeorient}
Q(X_1,X_1)=-1,\qquad Q(X_{\alpha},X_{\alpha})=1,\ \  \alpha=2,3,4, \qquad Q(X_i,X_j)=0\ \ \text{if}\ \ i\neq j.\end{equation}

Given a sub-Lorentzian metric $Q_x$ at $x\in\mathbf H$, we can define a linear mapping $g_x\colon
T^*\mathbf H_x\to T\mathbf H_x$ as follows: for $\xi\in T^*\mathbf H_x$, the linear mapping
$Y\to \xi(Y)$ for $Y\in T_h\mathbf H_x$ is represented uniquely by the identity
\begin{equation}\label{eq:2.1}
   Q_x(Y,g_x(\xi))=\xi(Y) \quad \mbox{for all}\;\; Y\in T_h\mathbf H_x.
\end{equation}

The properties of $Q_x$ implies that $g_x$ varies
smoothly in $x$, is symmetric, and nondegenerate. If $Q$ is of index 1, then index of $g$ is also 1.

Conversely, given a symmetric nondegenerate linear operator
$g_x\colon T^*\mathbf H_x\to T\mathbf H_x$ with the image $T_h\mathbf H_x$, there is a unique
nondegenerate quadratic form $Q_x$ satisfying \eqref{eq:2.1}.

{\it Example.} Let us consider the sub-Lorentzian metric $Q_x$ and co-metric $g_x$ on the quaternion $\mathbb H$-type group:
\begin{gather*}
    \{Q_x\}_{\alpha\beta}=\left(\begin{matrix}
           -1& 0 &0 &0 \\
            0& 1 &0 &0 \\
            0& 0 &1 &0 \\
            0 &0 &0 &1
\end{matrix}\right),
\end{gather*}
{\footnotesize{\begin{gather*}
    g_x^{\alpha\beta}=\left(\begin{matrix}
           -1& 0 &0 &0 &-\frac{1}{2}x_2 &\frac{1}{2}x_4 &\frac{1}{2}x_3\vspace{1mm}\\
            0& 1 &0 &0 &-\frac{1}{2}x_1 &-\frac{1}{2}x_3 &\frac{1}{2}x_4\vspace{1mm}\\
            0& 0 &1 &0 &\frac{1}{2}x_4 &\frac{1}{2}x_2 &\frac{1}{2}x_1\vspace{1mm}\\
            0 &0 &0 &1 &-\frac{1}{2}x_3 &\frac{1}{2}x_1 &-\frac{1}{2}x_2\vspace{1mm}\\
            -\frac{1}{2}x_2 &-\frac{1}{2}x_1 &\frac{1}{2}x_4 &-\frac{1}{2}x_3 &\frac{1}{4}(x_1^2-x_2^2+x_3^2+x_4^2) &\frac{1}{2}x_2x_4 &\frac{1}{2}x_2x_3\vspace{1mm}\\
            \frac{1}{2}x_4 &-\frac{1}{2}x_3 &\frac{1}{2}x_2 &\frac{1}{2}x_1 &\frac{1}{2}x_2x_4 &\frac{1}{4}(x_1^2+x_2^2+x_3^2-x_4^2) &-\frac{1}{2}x_3x_4\vspace{1mm}\\
            \frac{1}{2}x_3 &\frac{1}{2}x_4 &\frac{1}{2}x_1 &-\frac{1}{2}x_2 &\frac{1}{2}x_2x_3 &-\frac{1}{2}x_3x_4 &\frac{1}{4}(x_1^2+x_2^2-x_3^2+x_4^2)
      \end{matrix}\right).
\end{gather*}}}

The latter matrix is symmetric, thus, it can be diagonalised to
\begin{gather*}
    \widetilde{g}^{\alpha\beta}=\left(\begin{matrix}
           -1& 0 &0 &0 &0 &0 &0\\
            0& 1 &0 &0 &0 &0 &0\\
            0& 0 &1 &0 &0 &0 &0\\
            0 &0 &0 &1 &0 &0 &0\\
            0 &0 &0 &0 &0 &0 &0\\
            0 &0 &0 &0 &0 &0 &0\\
            0 &0 &0 &0 &0 &0 &0
      \end{matrix}\right).
\end{gather*} Therefore, the index of the co-metric $g^{\alpha\beta}$ is $1$ similarly to the index of the initial quadratic form $Q$.

Let us define the causal character of $\mathbf H$. Fix a point $p\in \textbf{H}$.
Denote by $\Omega_p$ the set of all horizontal curves $c(s)\colon [0,1]\to \textbf{H}$ starting from $c(0)=p$. A horizontal vector $v\in T_h\textbf{H}_p$ is called timelike if $Q(v,v)<0$, spacelike if $Q(v,v)>0$ or $v=0$, null if $Q(v,v)=0$ and $v\neq0$, nonspacelike if $Q(v,v)\leqslant 0$. A horizontal curve is called timelike if its tangent is timelike; similarly for spacelike, null and nonspacelike curves.

By the time-orientation of a sub-Lorentzian manifold $(\textbf{H}, T_h\textbf{H},Q)$ we mean a continuous horizontal timelike vector field on $\textbf{H}$. The time orientation of $(\textbf{H}, T_h\textbf{H},Q)$ is given by the vector $X_1$ according to~\eqref{timeorient}. Then a nonspacelike $v\in T_h\textbf{H}_p$ is called future directed if $Q(v,X_1(p))<0$, and it is called past directed if $Q(v,X_1(p))>0$. Throughout this paper f.d. stands for "future directed", t. for "timelike", and nspc. for "nonspacelike".

For an open set $U$ and fixed $p\in U$, we define two reachable sets: $I^+(p,U)$ $($resp. $J^+(p,U))$
is the set of all points $q\in U$ that can be reached from $p$ along a t.f.d. (resp. nspc.f.d.) curve contained in $U$. In the Lorentzian geometry $I^+(p,U)$  is called the chronological future of $p$ (with respect to $U$); similarly, $J^+(p,U)$ is called the causal future of $p$ (with respect to $U$).

We define a length of a nspc. curve $c\colon [\alpha,\beta]\to \textbf{H}$ in the following way
$$L(c)=\int\limits_{\alpha}^{\beta}|Q(\dot{c},\dot{c})|^{\frac{1}{2}}\,dt.$$

Let $\varphi\colon U\to \mathbb{R}$ be a smooth function on an open set $U$ in $\textbf{H}$. The horizontal gradient $\nabla_h\varphi$ of $\varphi$ is a smooth horizontal vector field on $U$ such that for each $p\in U$ and $v\in T_h\mathbf{H}_p$ we have $(\partial_v\varphi)(p)=Q(v,\nabla_H\varphi(p))$. Locally it can be written in the following way
$$\nabla_h\varphi=-(X_1\varphi)X_1+\sum\limits_{j=2}^4(X_j\varphi)X_j,$$ where $X_1,\ldots,X_4$ is an orthonormal frame of $T_h\textbf{H}$ defined on $U$ with $X_1$ timelike.

The following proposition can be found in \cite{Groch1}.
\begin{proposition}
    For any normal neighbourhood $U$ of a point $p$\\
    $(i)$ $J^{+}(p, U)$ is a closed subset in $U$;\\
    $(ii)$ $cl (I^+(p,U))=J^+(p,U)$, where $cl$ stands for the closure with respect to $U$.
\end{proposition}

\section{Reachable sets for Quaternion $\mathbb{H}$-Type Groups}

In this section we describe the reachable set for quaternion $\mathbb{H}$-type group. We start from the case $p=0$ and $I^+(0)=I^+(0,\mathbf{H})$.  Thereto let us consider the family of functions $\eta_{\alpha}=-x_1^2+x_2^2+x_3^2+x_4^2+\alpha(|z_{1}|+|z_{2}|
+|z_{3}|)$, $\,|\alpha|\leqslant \dfrac{4}{\sqrt{3}}$ that represents the homogeneous norm on~$\mathbf H$. We want to describe the reachable set in terms of values of the function $\eta_{\alpha}$. We present the calculations for the positive values of $z_{1},\,z_{2},\,z_{3}$, other cases considered analogously. The horizontal gradient of $\eta_{\alpha}$ is
\begin{gather}
    \nabla_{h}\eta_{\alpha} = \left[2x_1-\dfrac{1}{2}\alpha(x_2-x_3-x_4)\right]X_1
    + \left[2x_2+\dfrac{1}{2}\alpha(-x_1-x_3+x_4)\right]X_2\notag\\
    + \left[2x_3+\dfrac{1}{2}\alpha(x_1+x_2+x_4)\right]X_3
    + \left[2x_4+\dfrac{1}{2}\alpha(x_1-x_2-x_3)\right]X_4.\notag
\end{gather}

Consider $\eta_0$ and the set \begin{equation}\label{gamma0}\Gamma_{0}=\{\eta_{0}<0,\,x_1>0\}.\end{equation}
\begin{theorem}\label{teor1}
Let $I^+(0)$ denote the reachable set for $\mathbf H$ and $\Gamma_0$ be as in~\eqref{gamma0}. Then $$I^+(0)\subset\Gamma_0$$
\end{theorem}
\begin{proof}

Since $Q(\nabla_h\eta_0,X_1)=-2x_1$ and $Q(\nabla_h\eta_0,\nabla_h\eta_0)=4\eta_0$, we conclude that $\eta_0$ is t.f.d. in $\Gamma_0$. We claim that
\begin{equation}
\label{eq:11}J^+(0)\subset \overline{\Gamma}_0=\{\eta_0\leqslant 0,x_1\geqslant 0\},
\end{equation}
where we set $J^{+}(0)=J^+(0,\mathbf H)$. Let $p\in J^{+}(0)$. It means that there exists $c(t)$ nspc.~f.~d. $c\colon[0,T]\to \mathbb{R}^7$, $c(0)=0$, $c(T)=p$. Since $c$ is future directed it follows that $$Q\big(\dot{c}(t),X_1(c(t))\big)=-\dot{x}_1(c(t))<0.$$ Hence the function $x_1(c(t))$ increases from $x_1(0)=0$ along $c(t)$, and therefore, $x_1(c(T))=x_1(p)>0$. The nonspacelikeness  of $c$ implies that $$Q(\dot{c},\dot{c})=-\dot{x}_1^2+\dot{x}_2^2+\dot{x}_3^2+\dot{x}_4^2\leqslant 0.$$ We need to show that $x_1(T)\geqslant \sqrt{x_2^2(T)+x_3^2(T)+x_4^2(T)}$, taking into account that $x_1(0)=x_2(0)=x_3(0)=x_4(0)=0$. It will be shown if we prove the inequality
\begin{equation}\label{eq:111}
\sqrt{\dot{x}(t)_2^2+\dot{x}(t)_3^2+\dot{x}(t)_4^2}\geqslant \frac{d}{dt}\Big(\sqrt{x_2^2+x_3^2+x_4^2}\Big)=\dfrac{x_2\cdot\dot{x}_2+x_3\cdot\dot{x}_3+x_4\cdot\dot{x}_4}{\sqrt{x_2^2+x_3^2+x_4^2}},
\end{equation} since in this case $$\dot x_1(t)\geq \frac{d}{dt}\Big(\sqrt{x_2^2+x_3^2+x_4^2}\Big)\ \ \text{and}\ \ x_1(0)=\big(\sqrt{x_2^2+x_3^2+x_4^2}\big)(0)=0$$ $$\ \ \Longrightarrow\ \ x_1(T)\geq \big(\sqrt{x_2^2+x_3^2+x_4^2}\big)(T).$$
Consider the case $x_2\cdot\dot{x}_2+x_3\cdot\dot{x}_3+x_4\cdot\dot{x}_4\geqslant0$, because otherwise the inequality~\eqref{eq:111} is obvious. Squaring both parts of \eqref{eq:111} and simplifying, we get that~\eqref{eq:111} holds if and only if $$(\dot{x}_2x_3-\dot{x}_3x_2)^2+(\dot{x}_2x_4-\dot{x}_4x_2)^2+(\dot{x}_3x_4-\dot{x}_4x_3)^2\geqslant0.$$ Statement \eqref{eq:11} is proved.

Now we look at the behaviour of any nspc.f.d. curve which projects onto the set $$\left\{ x_1= \sqrt{x_2^2+x_3^2+x_4^2},\,x_1>0\right\}.$$ Due to the chain of inequalities $$\dot{x}_1=\frac{d}{dt}\Big(\sqrt{x_2^2+x_3^2+x_4^2}\Big) \leq
\sqrt{\dot{x}_2^2+\dot{x}_3^2+\dot{x}_4^2}\leq\dot{x}_1,$$ such curves must be null curves. In addition, since $\nabla_h\eta_0$ is nspc.f.d. on $\overline{\Gamma}_0$, we see that for any t.f.d. curve $c\colon [0,T]\to \mathbb{R}^7$, $c(0)=0$, the derivative along the curve $c(t)$ $$\dfrac{d}{dt}\eta_0(c(t))=Q(\nabla_h\eta_0(c(t)),\dot c(t))<0,$$ and the function $t\longrightarrow \eta_0(c(t))$ is decreasing a.~e. It follows that $I^+(0)\subset \Gamma_0$.
\end{proof}

\begin{corollary} Following the notations of the previous theorem we have
$$\Gamma_0\cap\{z_1=0,z_2=0,z_3=0\}= I^+(0)\cap\{z_1=0,z_2=0,z_3=0\}.$$
\end{corollary}

\begin{proof}
Consider the straight lines
\begin{equation}
\label{eq:2.3}
    \begin{cases}
        x_1=t\cosh \varphi,\\
        x_2=t\sinh\varphi\sin\psi\cos\theta,\\
        x_3=t\sinh\varphi\sin\psi\sin\theta,\\
        x_4=t\sinh\varphi\cos\psi, \quad t>0,
    \end{cases}
\end{equation}
where $\psi\in[0,2\pi]$, $\theta\in[0,\pi]$, $\varphi\in[-\infty,+\infty]$ are constants. Any of these straight lines is a t.f.d. curve. They fill up the interior of $\Gamma_0\cap\{z_1=0,z_2=0,z_3=0\}$. So, we conclude that $$\Gamma_0\cap\{z_1=0,z_2=0,z_3=0\}\subset I^+(0)\cap\{z_1=0,z_2=0,z_3=0\}.$$ The inverse inclusion follows from Theorem~\ref{teor1}.
\end{proof}

Let us calculate sets, where $\nabla_h\eta_{\alpha}$ is nspc.f.d. We have
$$Q(\nabla_h\eta_{\alpha}, X_1)=2(-x_1+\frac{\alpha}{4}(x_2-x_3-x_4))=2(-x_1+\frac{\sqrt 3\alpha}{4}(\overrightarrow{y}\cdot\overrightarrow{n})),$$ \begin{eqnarray*}
Q(\nabla_h\eta_{\alpha},\nabla_h\eta_{\alpha}) & = &
4\Big((-1+\frac{3\alpha^2}{16})(x_1^2-x_2^2-x_3^2-x_4^2)\Big)\\ & + & \frac{\alpha^2}{2}\big((x_2+x_3)^2+(x_2+x_4)^2+(x_3-x_4)^2\big)\\ & = &
4\Big((-1+\frac{3\alpha^2}{16})(x_1^2-|y|^2)+\frac{3\alpha^2}{8}|y\times n|^2\Big),
\end{eqnarray*}
where we used the notation $\overrightarrow{y}=(x_2,x_3,x_4)$, $\overrightarrow{n}=\frac{1}{\sqrt 3}(1,-1,-1)$ and $``\times"$ is the vector product in $\mathbb R^3$ and $|\cdot|$ is the length of the vector with respect to the inner product in $\mathbb R^3$.
We see that $\nabla_h\eta_{\alpha}$ is f.d for $|\alpha|<\frac{4}{\sqrt 3}$ in $\Gamma_0$. Indeed, if we write the condition of future directness as $$-x_1+\frac{\alpha}{4}(x_2-x_3-x_4)=\langle \overrightarrow{N},\overrightarrow{x}\rangle<0, \quad \overrightarrow{N}=(1,\frac{\alpha}{4},-\frac{\alpha}{4},-\frac{\alpha}{4}),\,\,\overrightarrow{x}=(x_1,x_2,x_3,x_4),$$ then it holds if and only if $$\langle \overrightarrow{N},\overrightarrow{N}\rangle<0,\quad\langle \overrightarrow{x},\overrightarrow{x}\rangle<0.$$ The latter two conditions satisfy in $\Gamma_0$ for $|\alpha|<\frac{4}{\sqrt 3}$.

We introduce the notations $$\Gamma_{\alpha}=\{\eta_{\alpha}<0,\ \ x_1>0\}$$ and $$A_{\alpha}=\{(-1+\frac{3\alpha^2}{16})(x_1^2-|\overrightarrow{y}|^2)+\frac{3\alpha^2}{8}|\overrightarrow{y}\times \overrightarrow{n}|^2<0,\quad x_1>0\}.$$ If $\alpha=0$, then $A_{0}=\Gamma_0$, otherwise $A_{\alpha}\subset\Gamma_0$. If $\alpha=\frac{4}{\sqrt 3}$ then $A_{\frac{4}{\sqrt 3}}=\{\emptyset\}$ and $\overline{A}_{\frac{4}{\sqrt 3}}=\{x_2=-x_3=-x_4,\ \ x_1\geq 0\}$.

\begin{remark}
The set $A_{\alpha}$ has the following geometric sense. Let $$Q=\{Q_{ij}\}=\begin{pmatrix}
 -1  & 0 & 0 & 0\\
0  & 1   & 0 & 0\\
0  & 0  & 1 & 0\\
0  & 0  & 0 & 1
\end{pmatrix},$$ let $$V=\{V_{kl}\}=\begin{pmatrix}
 \frac{1}{4}\big(x_1^2-x_2^2+x_3^2+x_4^2\big)  & \frac{1}{2}x_2x_4   & \frac{1}{2}x_2x_3 \\
\frac{1}{2}x_2x_4  & \frac{1}{4}\big(x_1^2+x_2^2+x_3^2-x_4^2\big)   & -\frac{1}{2}x_3x_4 \\
\frac{1}{2}x_2x_3  & -\frac{1}{2}x_3x_4  & \frac{1}{4}\big(x_1^2+x_2^2-x_3^2+x_4^2\big)
\end{pmatrix}$$ be quadratic forms and let $\overrightarrow{x}=(x_1,x_2,x_3,x_4)$, $\overrightarrow{z}=(z_1,z_2,z_3)$ be  vectors. Then $A_{\alpha}$ can be written as $$A_{\alpha}=\{Q(\overrightarrow{x},\overrightarrow{x})+V(\overrightarrow{z},\overrightarrow{z})<0,\ \ x_1>0\},\ \ \text{with}\ \ z_1=z_2=z_3=\alpha.$$ We notice that both of the quadratic forms appeared in the co-metric $g^{\alpha\beta}$.
\end{remark}

\begin{theorem}
We have $$I^+(0)\cap A_{\alpha}\subset\Gamma_{\alpha}.$$
\end{theorem}
\begin{proof}
Let $p\in I^+(0)\cap A_{\alpha}$. Then there is a timelike f.d. curve $c:(0,T)\to \Gamma_0\cap A_{\alpha}$, $c(T)=p$. Since the vector $\nabla_h\eta_{\alpha}$ is timelike and f.d. in the same set we conclude that the derivative $\frac{d}{dt}\eta_{\alpha}(c(t))=Q(\nabla_h\eta_{\alpha},\dot c)$ is negative and the function $\eta_{\alpha}(c(t))$ decrease along $c(t)$, that yields $\eta_{\alpha}(p)<0$ jointly with $x_1(p)>0.$
\end{proof}

The complete description of the reachable set $I^+(0)$ is complicated, therefore we study its subsets.

\begin{proposition}
Let $B$ be one of the following sets $\{x_3=x_4=0\}$, $\{x_2=x_3=0\}$, $\{x_2=x_4=0\}$. Then $$I^+(0)\cap B=\Gamma_4\cap B.$$
\end{proposition}
\begin{proof}
If $x_3=x_4=0$, then for any horizontal curve the horizontality conditions
\begin{eqnarray}\label{eq:horizontality}
\dot z_1 & = & \frac{1}{2}(+x_2 \dot x_1-x_1\dot x_2+x_4\dot x_3-x_3\dot x_4),\notag\\
\dot z_2 &= & \frac{1}{2}(-x_4 \dot x_1-x_3\dot x_2+x_2\dot x_3+x_1\dot x_4),\\
\dot z_3 & = & \frac{1}{2}(-x_3 \dot x_1+x_4\dot x_2+x_1\dot x_3-x_2\dot x_4)\notag
\end{eqnarray} implies that $\dot z_1=\frac{1}{2}(+x_2 \dot x_1-x_1\dot x_2)$ and $\dot z_2=\dot z_3=0$.
Since $z_1(0)=z_2(0)=z_3(0)=0$, we conclude that this case is reduced to the Heisenberg sub-Lorentzian manifold and we can apply results of~\cite{Groch1}. Other cases are obtained analogously.
\end{proof}

\section{Hamiltonian and geodesics in Heisenberg group with Lorentzian metric}

In this section we study geodesics on the one dimensional
Heisenberg group equipped with the Lorentzian metric. We remind
that the Heisenberg group $\mathbb H^1$ is the space $\mathbb R^3$
furnished with the non-commutative law of multiplication
$$(x,y,z)(x^{\prime},y^{\prime},z^{\prime})=\big(x+x^{\prime},y+y^{\prime},z+z^{\prime}+\frac{1}{2}(yx^{\prime}-xy^{\prime})\big).$$
The two dimensional horizontal bundle $T_h\mathbb H^1$ is given as
a span of left invariant vector fields
$$X=\frac{\partial}{\partial
x}+\frac{1}{2}y\frac{\partial}{\partial z},\qquad
Y=\frac{\partial}{\partial y}-\frac{1}{2}x\frac{\partial}{\partial
z},\qquad [X,Y]=Z=\frac{\partial}{\partial z}.$$ We suppose that the
Lorentzian metric $Q$ is defined on $T_h\mathbb H^1$ by
$$Q(X,X)=-1,\quad Q(Y,Y)=1,\quad Q(X,Y)=0.$$ The time orientation
is given by the horizontal vector field $X$. Thus the triple
$(\mathbb R^3, T_h\mathbb H^1, Q)$ is called the Heisenberg group
with the Lorentzian metric, and to differ it from the classical case we
use the notation $\mathbb H^1_L$. The reachable set for the Heisenberg
group with the Lorentzian metric was studied in~\cite{Groch1, Groch2}.

To study geodesics we apply the Hamiltonian method. The
Hamiltonian is defined as a symbol of the corresponding wave
equation formed by the left invariant horizontal vector fields.
The definition of geodesics in sub-Riemannian geometry differs
from the definition of geodesics in the Riemannian geometry.
Geodesics in the sub-Riemannian geometry (and we shall use the
same definition) are projections of solutions of the corresponding
Hamiltonian system onto the underlying manifold. The interesting
feature of the sub-Riemannian geometry that even locally there is
no uniqueness of geodesics. The study of number of geodesics on
Heisenberg group with positively definite metric can be found, for
instance, in~\cite{CChGr3,Gaveau}. We apply the method of
Hamiltonian mechanics in order to calculate geodesics of different
causal characters and to study the uniqueness problem.

The governed operator is an analogue of the sub-Laplacian which
has the form
\begin{gather}
  2\bigtriangleup_H=-X^2+Y^2=-\left(\dfrac{\partial }{\partial x}+\dfrac{1}{2}y\dfrac{\partial}{\partial z}\right)^2+\left(\dfrac{\partial}{\partial   y}-\dfrac{1}{2}x\dfrac{\partial}{\partial z}\right)^2\notag\\
   =-\dfrac{\partial^2 }{\partial x^2}-\dfrac{1}{4}y^2\dfrac{\partial^2 }{\partial z^2}
   -y\dfrac{\partial }{\partial x}\dfrac{\partial }{\partial z} +\dfrac{\partial^2 }{\partial y^2}+\dfrac{1}{4}x^2\dfrac{\partial^2 }{\partial z^2} -x\dfrac{\partial }{\partial y} \dfrac{\partial }{\partial z}\notag\\
   =\left(-\dfrac{\partial^2 }{\partial x^2}+\dfrac{\partial^2 }{\partial y^2}\right)-\dfrac{1}{4}(-x^2+y^2)\dfrac{\partial^2 }{\partial z^2}
   -y\dfrac{\partial }{\partial x} \dfrac{\partial }{\partial z}-x\dfrac{\partial }{\partial y} \dfrac{\partial }{\partial z}.\notag
\end{gather}
Then the associated Hamiltonian function $H(\xi,\eta,\theta,x,y,z)$ becomes
$$H=\dfrac{1}{2}(-\xi^2+\eta^2)-\dfrac{1}{8}(-x^2+y^2)\theta^2+\dfrac{1}{2}(-y\xi-x\eta)\theta,$$
where we use the notations $\xi=\frac{\partial}{\partial x}$,
$\eta=\frac{\partial}{\partial y}$,
$\theta=\frac{\partial}{\partial z}$. The corresponding
Hamiltonian system is
\begin{equation}\label{hamheis}
\begin{cases}
    &\dot{x}=\dfrac{\partial H}{\partial\xi}=-\xi-\dfrac{y\theta}{2},\\[3mm]
    &\dot{y}=\dfrac{\partial H}{\partial\eta}=\eta-\dfrac{x\theta}{2},\\[3mm]
    &\dot{z}=\dfrac{\partial H}{\partial\theta}=-\dfrac{1}{4}(-x^2+y^2)\theta
    +\dfrac{1}{2}(-y\xi-x\eta),\\[3mm]
    &\dot{\xi}=-\dfrac{\partial H}{\partial x}=-\dfrac{1}{4}x\theta^2+\dfrac{1}{2}\eta\theta,\\[3mm]
    &\dot{\eta}=-\dfrac{\partial H}{\partial y}=\dfrac{1}{4}y\theta^2+\dfrac{1}{2}\xi\theta,\\[3mm]
    &\dot{\theta}=-\dfrac{\partial H}{\partial z}=0. \end{cases} \end{equation}
The initial data is $x(0)=y(0)=z(0)=0$, $\xi(0)=\xi_0$,
$\eta(0)=\eta_0$, $\theta(0)=\theta$. We are interested in finding
the projection of the solution of~\eqref{hamheis} onto the
$(x,y,z)$-space. We reduce the system~\eqref{hamheis} to the
system containing only $(x,y,z)$ coordinates withthe initial
data zero. If we express $\xi$ and $\eta$ from the first two equations
and substitute them in the expressions for $\dot{\xi}$ and
$\dot{\eta}$, we obtain $\dot\xi=\frac{1}{2}\dot y\theta$,
$\dot\eta=-\frac{1}{2}\dot x\theta$. Taking into account that
$\theta$ is constant, we differentiate the first two equations and
replace $\dot\xi$ and $\dot\eta$ there. We get
\begin{gather}
   \begin{cases}
      \ddot{x}=-\dot{y}\theta,\\
      \ddot{y}=-\dot{x}\theta
   \end{cases}\notag
\end{gather}
or
\begin{gather}
   \left(%
\begin{array}{c}
  \ddot{x} \\
  \ddot{y} \\
\end{array}%
\right)=\left(%
\begin{array}{cc}
  0 & -\theta \\
  -\theta & 0 \\
\end{array}%
\right)\left(%
\begin{array}{c}
  \dot{x} \\
  \dot{y} \\
\end{array}%
\right).\notag
\end{gather}
We are looking for the solution $x=x(t)$, $y=y(t)$,
$t\in[-\infty,+\infty]$, satisfying $x(0)=0$, $y(0)=0$ and
$\dot{x} (0)=\dot{x}_0=-\xi_0$, $\dot{y}(0)=\dot{y}_0=\eta_0$.
Eigenvalues and eigenvectors for the matrix
$M:=\left(%
\begin{array}{cc}
  0 & \theta \\
  \theta & 0 \\
\end{array}%
\right)$ are $\lambda=\pm|\theta|$ and $\left(%
\begin{array}{c}
  1 \\
  -1 \\
\end{array}%
\right)$, $\left(%
\begin{array}{c}
  1 \\
  1 \\
\end{array}%
\right)$ respectively. Therefore,
\begin{gather}
   \begin{cases}
       \dot{x}(t)=c_1e^{|\theta| t}+c_2e^{-|\theta| t},\quad\dot{x}_0=c_1+c_2,\\
       \dot{y}(t)=-c_1e^{|\theta| t}+c_2e^{-|\theta|
       t},\quad\dot{y}_0=-c_1+c_2.
   \end{cases}\notag
\end{gather}
It follows that
\begin{equation}\label{xyderivative}
   \begin{cases}
       \dot{x}(t)=\dot{x}_0\cosh(|\theta| t)-\dot{y}_0\sinh(|\theta| t),\\
       \dot{y}(t)=-\dot{x}_0\sinh(|\theta| t)+\dot{y}_0\cosh(|\theta| t).
   \end{cases}
\end{equation}
Finally, we obtain
\begin{equation}\label{xycoordinate}
\begin{cases}
    x(t)=\dfrac{\dot{x}_0}{|\theta|}\sinh(|\theta| t)-\dfrac{\dot{y}_0}{|\theta|}\big(\cosh(|\theta| t)-1\big),\\
    y(t)=-\dfrac{\dot{x}_0}{|\theta|}\big(\cosh(|\theta| t)-1\big)+\dfrac{\dot{y}_0}{|\theta|}\sinh(|\theta| t). \end{cases}
\end{equation}

\begin{lemma}
 If a geodesic is timelike future directed (past directed) at $t=0$, then it remains timelike and future directed (past directed) for all $t\in[0,\pm\infty]$.
\end{lemma}
\begin{proof}
Since geodesics satisfy the Hamiltonian system~\eqref{hamheis},
the speed is preserved along geodesics and this implies that the
causality character does not change. This also follows
from~\eqref{xyderivative} because $-\dot x^2(t)+\dot
y^2(t)=-\dot x^2_0+\dot y^2_0$. We remind that a geodesic is
future directed if $\dot x(t)>0$. Let us show that if we start
from $\dot x_0>0$ and $-\dot x_0<\dot y_0<\dot x_0$, then the
geodesic remains future directed for all $t\in (0,\pm\infty)$. If
$0<\dot y_0<\dot x_0$, then we conclude that $\dot x=\dot x_0\cosh
(|\theta|t)-\dot y_0\sinh (|\theta|t)>0$ for $t>0$ because of the inequality
$\sinh (|\theta|t)<\cosh (|\theta|t)$. If $-\dot x_0<\dot y_0<0$
then we take $t\in(0,-\infty)$ and get future directed geodesics
because of $-\sinh (|\theta|t)<\cosh (|\theta|t)$. We conclude
that for $\dot x_0>0$, and $|\dot y_0|<\dot x_0$ we obtain future
directed geodesics for $t\in (0,\pm\infty)$.
\end{proof}

Using the condition of horizontality, we get
\begin{equation*}
\dot{z}  = \dfrac{1}{2}(y\dot{x}-x\dot{y})=
\dfrac{1}{2|\theta|}(-\dot{x}_0^2+\dot{y}_0^2)(1-\cosh (|\theta| t)).
\end{equation*} Notice that the derivative $\dot z$  is positive for timelike geodesics, for lightlike curves we have $\dot z=0$, and spacelike geodesics satisfy $\dot z<0$ for any values of $t$.
Integrating, we get
\begin{gather}\label{zgeneral}
    z(t)=\dfrac{\|v_0\|^2}{2\theta^2}\left(|\theta|t-\sinh (|\theta| t)\right),\qquad \|v_0\|^2=-\dot{x}_0^2+\dot{y}_0^2.
\end{gather} The value of $z$-coordinate for timelike future directed geodesics starting from the origin is positive for positive value of $t$ and negative for negative value of $t$, lightlike curves has vanishing $z$-coordinate and spacelike curves admits the negative values for $t>0$ and positive for $t<0$.

\begin{figure}[ht]
\centering \scalebox{0.8}{\includegraphics{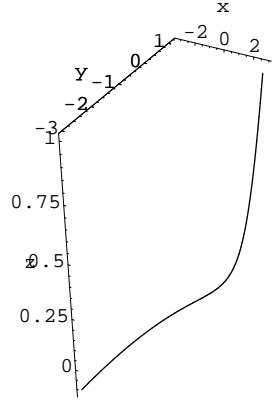}}
\caption[]{The graph of timelike geodesic passing through the
origin}\label{fig:timelike}
\end{figure}

Further, introducing the notation $\|h\|^2=-x^2+y^2$ we get
\begin{equation}\label{horcoord}\|h\|^2=-x^2+y^2=\frac{4\|v_0\|^2}{\theta^2}\sinh^2\big(\frac{|\theta|t}{2}\big).
\end{equation}
If a geodesic is timelike, then its projection onto
$(x,y)$-plane lies inside the domain $-x^2+y^2<0$. Lightlike
geodesics are projected into the set $-x^2+y^2=0$ and spacelike
geodesics are projected into the domain $-x^2+y^2>0$, see
Figure~\ref{fig:timelike}.

\begin{lemma} Let $\gamma(t)=\big(x(t),y(t),z(t)\big)$ be a
geodesic. Then
\begin{equation*}
-x^2(t)+y^2(t)+4|z(t)|\qquad\begin{cases}
    <0\quad\text{if}\quad\gamma(t)\quad\text{is timelike} ,\\
    =0\quad\text{if}\quad\gamma(t)\quad\text{is lightlike} ,\\
    >0\quad\text{if}\quad\gamma(t)\quad\text{is spacelike} . \end{cases}
\end{equation*}
\end{lemma}
\begin{proof}
We use the notation $\|h(t)\|^2=-x^2(t)+y^2(t)$. If $z>0$, then the value
$$\|h(t)\|^2+4 |z(t)|=\|h(t)\|^2+4 z(t)=\frac{2\|v_0\|^2}{\theta^2}(-1+|\theta|t+e^{-|\theta|t})$$
is negative for timelike curves and positive for spacelike curves,
because $-1+|\theta|t+e^{-|\theta|t}>0$ for all $t\neq 0$. The value of 
$\|h\|^2+4 z$ vanishes for lightlike curves. If $z<0$, then the
expression $$\|h(t)\|^2+4 |z(t)|=\|h(t)\|^2-4
z(t)=\frac{2\|v_0\|^2}{\theta^2}(-1-|\theta|t+e^{|\theta|t})$$ is
still negative for timelike curves and is positive for spacelike
geodesics, because $-1-|\theta|t+e^{|\theta|t}>0$ for all $t\neq 0$.
The value $\|h\|^2-4 z$ vanishes for lightlike curves.
\end{proof}

\begin{lemma}
The equation
\begin{equation}\label{eq:12}\frac{4z}{-x^2+y^2}=\frac{\tau}{\sinh^2(\tau)}-\coth(\tau)\end{equation}
has a unique solution $\tau$ for given $(x,y,z)$  if
$-1<\frac{4z}{-x^2+y^2}<1$, and has no solution otherwise.
\end{lemma}
\begin{proof}
The function $\mu(\tau)=\frac{\tau}{\sinh^2(\tau)}-\coth(\tau)$ is
strictly decreasing in the interval $(-\infty,+\infty)$ from $1$
to $-1$, see Figure~\ref{fig:graph1}. It proves the lemma.

\begin{figure}[ht]
\centering \scalebox{0.7}{\includegraphics{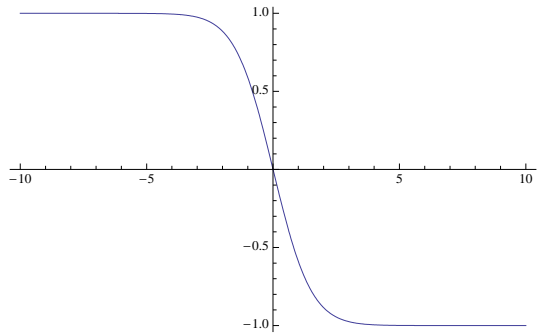}}
\caption[]{The graph of function $\mu(\tau)$}\label{fig:graph1}
\end{figure}

\end{proof}

We would like to draw the reader's attention to the function $\mu$. It is
intriguing analogue of the function used by Gaveau~\cite{Gaveau}
to study geodesics and geometry in general on the classical
Heisenberg group. The classical counterpart $\widetilde \mu$ of
the function $\mu$ has the form $\widetilde
\mu(\tau)=\frac{\tau}{\sin^2\tau}-\cot\tau$, and the classical
Heisenberg analogue of the equation~\eqref{eq:12},
$$\frac{4|z|}{x^2+y^2}=\frac{\tau}{\sin^2\tau}-\cot\tau$$ has more
then 1 solution.

We state the following theorems describing the reachable set by
geodesics starting from the origin.

\begin{theorem}\label{th:timelike}
Let $A=(x,y,z)$ be a point such that $x>0$ $(x<0)$,$-x^2+y^2<0$,
$\frac{4|z|}{x^2-y^2}<1$. Then there is a unique timelike future
directed $($past directed$)$ geodesic, joining $O=(0,0,0)$ with
the point $A$. Let $\theta$ be a solution of the equation
\begin{equation}\label{mu}\frac{4z}{-x^2+y^2}=\frac{|\theta|/2}{\sinh^2(|\theta|/2)}-\coth(|\theta|/2).
\end{equation}
Then the equations of timelike future directed geodesic $\gamma: [0,1]\to \mathbb H^1_L$ are
\begin{equation}\label{xtime}
\begin{cases}
    x(t)=\sinh^2(\frac{|\theta|}{2}t)\Big(x\big(\coth(\frac{|\theta|}{2}t)\coth(\frac{|\theta|}{2})-1\big)+y\big(\coth(\frac{|\theta|}{2}t)-\coth(\frac{|\theta|}{2})\big)\Big),\\
    y(t)=\sinh^2(\frac{|\theta|}{2}t)\Big(y\big(\coth(\frac{|\theta|}{2}t)\coth(\frac{|\theta|}{2})-1\big)+x\big(\coth(\frac{|\theta|}{2}t)-\coth(\frac{|\theta|}{2})\big)\Big). \end{cases}
\end{equation}

\begin{equation}\label{ztime}z(t)=z\frac{|\theta|t-\sinh(|\theta|t)}{|\theta|-\sinh(|\theta|)}\end{equation}
The square of the length of these geodesics is $$l^2=\theta^2\frac{|\|h(1)\|^2|+4|z|}{2\big(||\theta|-\sinh(|\theta|)|+2\sinh^2(|\theta|/2)\big)}.$$
\end{theorem}

\begin{figure}[ht]
\centering \scalebox{0.7}{\includegraphics{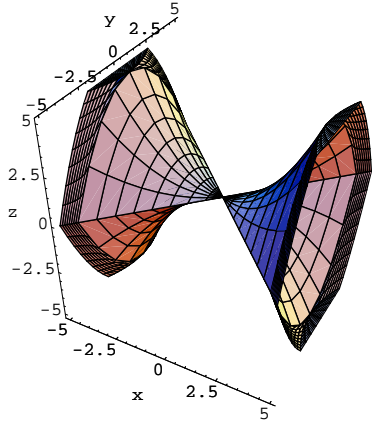}}
\caption[]{The domain reachable from the origin by timelike
geodesics}\label{timelikeregion}
\end{figure}

\begin{proof}
Let us suppose that geodesics joining $O=(0,0,0)$ with $A$ are
parametrized on the interval $[0,1]$. We fix the solution of the
equation~\eqref{mu} for a given point $A=(x,y,z)$, $x>0$,
$-1<\frac{4z}{x^2-y^2}<1$. Put $t=1$ in~\eqref{zgeneral} and find
$\frac{\|v_0\|^2}{2\theta^2}=\frac{z}{|\theta|-\sinh|\theta|}$.
Then we substitute the last expression in~\eqref{zgeneral} and
get~\eqref{ztime}.

To obtain~\eqref{xtime} we do essentially the same. Set $t=1$ in the expressions for $x(t)$ and $y(t)$ to find $\dot x_0$, $\dot y_0$. We get $$\dot x_0=\frac{|\theta|}{2}(x\coth(\frac{|\theta|}{2})+y),\quad \dot y_0=\frac{|\theta|}{2}(y\coth(\frac{|\theta|}{2})+x).$$ Substituting these expressions in~\eqref{xycoordinate}, we get~\eqref{xtime}.

Let us calculate the length of timelike geodesics. Expressing $\|v_0\|^2$ from~\eqref{horcoord}, we get $$\|v_0\|^2=\frac{\|h(1)\|^2\theta^2}{4\sinh^2(\frac{|\theta|}{2})}.$$ Then the length of geodesics is
\begin{equation}\label{eq:length}
l=\int_0^1\sqrt{|\|v_0\|^2|}\,dt=\frac{\sqrt{|\|h(1)\|^2}||\theta|}{2\sinh(\frac{|\theta|}{2})}\,dt.
\end{equation}
We calculate, making use of equation~\eqref{mu}
$$|\|h(1)\|^2|+4|z|=|\|h(1)\|^2|\Big(\Big|\frac{|\theta|/2}{\sinh^2(|\theta|/2)}-\coth(|\theta|/2)\Big|+1\Big).$$
Substituting the expression for the length, we get
$$l^2=\theta^2\frac{|\|h(1)\|^2|+4|z|}{2\big(||\theta|-\sinh(|\theta|)|+2\sinh^2(|\theta|/2)\big)}.$$
In the Figure~\ref{timelikeregion} we present the domain reachable
from the origin by timelike geodesics.
\end{proof}

In the case of spacelike curve we prove the analogous theorem
\begin{theorem}
Let $B=(x,y,z)$ be a point such that $-x^2+y^2>0$, $-1<\frac{4z}{-x^2+y^2}<1$, then there is a unique spacelike geodesic, joining $O=(0,0,0)$ with the point $B$. Let $\theta$ satisfy the equation~\eqref{mu}
Then the equations of spacelike geodesics are given by~\eqref{xtime} and ~\eqref{ztime}
The square of the length of these geodesics is
$$l^2=\theta^2\frac{\|h(1)\|^2+4|z|}{2\big(||\theta|-\sinh(|\theta|)|+2\sinh^2(|\theta|/2)\big)}.$$
\end{theorem}

\begin{remark}
Notice that~\eqref{zgeneral} and~\eqref{horcoord} imply that the
horizontal lightlike geodesics starting from the origin satisfy the
equations $x(t)=\pm y(t)$, $z(t)=0$.
\end{remark}

\begin{remark}
The result of Theorem~\ref{th:timelike} gives the description of
the set reachable from the origin by timelike geodesics.
Grochowski~\cite{Groch1} proved that the region
$\frac{4|z|}{x^2-y^2}<1$, $x>0$ is reachable by horizontal
timelike future directed curves. Theorem~\ref{th:timelike} states
that this set is reachable by geodesics of the same causal
character.
\end{remark}

\section{Hamiltonian and geodesics for $\mathbf H$ with the Lorentzian metric}

In the present section we find the parametric equations of geodesics for quaternion $\mathbb H$-type groups furnished with the Lorentzian metric. The horizontal wave operator is
\begin{gather}
   2\triangle_H=-X_1^2+X_2^2+X_3^2+X_4^2\notag\\
    =-\dfrac{\partial^2}{\partial x_1^2}+\dfrac{\partial^2}{\partial
   x_2^2}+\dfrac{\partial^2}{\partial x_3^2}+\dfrac{\partial^2}{\partial x_4^2}
   +x_2x_4\dfrac{\partial^2}{\partial z_1\partial z_2}+x_2x_3\dfrac{\partial^2}{\partial z_1\partial
   z_3}-x_3x_4\dfrac{\partial^2}{\partial z_2\partial z_3}\notag\\+\dfrac{1}{4}\dfrac{\partial^2}{\partial
   z_1^2}(x_1^2-x_2^2+x_3^2+x_4^2)+\dfrac{1}{4}\dfrac{\partial^2}{\partial
   z_2^2}(x_1^2+x_2^2+x_3^2-x_4^2)+\dfrac{1}{4}\dfrac{\partial^2}{\partial
   z_3^2}(x_1^2+x_2^2-x_3^2+x_4^2)\notag\\
   +\dfrac{\partial}{\partial z_1}\left(-x_2\dfrac{\partial}{\partial x_1}-x_1\dfrac{\partial}{\partial x_2}
   +x_4\dfrac{\partial}{\partial x_3} -x_3\dfrac{\partial}{\partial x_4}\right)
   +\dfrac{\partial}{\partial z_2}\left(x_4\dfrac{\partial}{\partial x_1}-x_3\dfrac{\partial}{\partial x_2}
   +x_2\dfrac{\partial}{\partial x_3} +x_1\dfrac{\partial}{\partial
   x_4}\right)\notag\\
   +\dfrac{\partial}{\partial z_3}\left(x_3\dfrac{\partial}{\partial x_1}+x_4\dfrac{\partial}{\partial x_2}
   +x_1\dfrac{\partial}{\partial x_3} -x_2\dfrac{\partial}{\partial
   x_4}\right).\notag
\end{gather}

Then the associated Hamiltonian function $H(\xi,\theta,x,z)$ has
the following form
\begin{gather}
   H=\frac{1}{2}(-\xi_1^2+\xi_2^2+\xi_3^2+\xi_4^2)+\frac{1}{2}(x_2x_4\theta_1\theta_2+x_2x_3\theta_1\theta_3-x_3x_4\theta_2\theta_3)\notag\\
   +\dfrac{1}{8}\theta_1^2(x_1^2-x_2^2+x_3^2+x_4^2)+\dfrac{1}{8}\theta_2^2(x_1^2+x_2^2+x_3^2-x_4^2)
   +\dfrac{1}{8}\theta_3^2(x_1^2+x_2^2-x_3^2+x_4^2)\notag\\
   +\frac{1}{2}\theta_1(-x_2\xi_1-x_1\xi_2+x_4\xi_3-x_3\xi_4)+\frac{1}{2}\theta_2(x_4\xi_1-x_3\xi_2+x_2\xi_3+x_1\xi_4)\\ +\frac{1}{2}\theta_3(x_3\xi_1+x_4\xi_2+x_1\xi_3-x_2\xi_4).\notag
\end{gather}
The corresponding Hamiltonian system is
\begin{equation}
\begin{array}{l}
\vspace{1mm}
  \dot{x}_1=\dfrac{\partial H}{\partial \xi_1}=-\xi_1-\frac{1}{2}x_2\theta_1+\frac{1}{2}x_4\theta_2+\frac{1}{2}x_3\theta_3, \\ \vspace{1mm}
  \dot{x}_2=\dfrac{\partial H}{\partial \xi_2}=\xi_2-\frac{1}{2}x_1\theta_1-\frac{1}{2}x_3\theta_2+\frac{1}{2}x_4\theta_3, \\ \vspace{1mm}
  \dot{x}_3=\dfrac{\partial H}{\partial \xi_3}=\xi_3+\frac{1}{2}x_4\theta_1+\frac{1}{2}x_2\theta_2+\frac{1}{2}x_1\theta_3, \\ \vspace{1mm}
  \dot{x}_4=\dfrac{\partial H}{\partial \xi_4}=\xi_4-\frac{1}{2}x_3\theta_1+\frac{1}{2}x_1\theta_2-\frac{1}{2}x_2\theta_3,\\ \vspace{1mm}
  \dot{z}_1=\dfrac{\partial H}{\partial \theta_1}=\frac{1}{2}(x_2x_4\theta_2+x_2x_3\theta_3)+\frac{1}{8}\theta_1(x_1^2-x_2^2+x_3^2+x_4^2)
  +\frac{1}{2}(-x_2\xi_1-x_1\xi_2+x_4\xi_3-x_3\xi_4), \\ \vspace{1mm}
  \dot{z}_2=\dfrac{\partial H}{\partial \theta_2}=\frac{1}{2}(x_2x_4\theta_1-x_3x_4\theta_3)+\frac{1}{8}\theta_2(x_1^2+x_2^2+x_3^2-x_4^2)
  +\frac{1}{2}(x_4\xi_1-x_3\xi_2+x_2\xi_3+x_1\xi_4), \\ \vspace{1mm}
  \dot{z}_3=\dfrac{\partial H}{\partial \theta_3}=\frac{1}{2}(x_2x_3\theta_1-x_3x_4\theta_2)+\frac{1}{8}\theta_3(x_1^2+x_2^2-x_3^2+x_4^2)
  +\frac{1}{2}(x_3\xi_1+x_4\xi_2+x_1\xi_3-x_2\xi_4),\\ \vspace{1mm}
 \dot{\xi}_1=-\dfrac{\partial H}{\partial x_1}=-\frac{1}{4}x_1(\theta_1^2+\theta_2^2+\theta_3^2)+\frac{1}{2}\xi_2\theta_1-\frac{1}{2}\xi_4\theta_2-\frac{1}{2}\xi_3\theta_3, \\ \vspace{1mm}
  \dot{\xi}_2=-\dfrac{\partial H}{\partial x_2}=-\frac{1}{4}x_2(-\theta_1^2+\theta_2^2+\theta_3^2)+\frac{1}{2}\xi_1\theta_1-\frac{1}{2}\xi_3\theta_2+\frac{1}{2}\xi_4\theta_3
  -\frac{1}{2}x_4\theta_1\theta_2-\frac{1}{2}x_3\theta_1\theta_3,  \\ \vspace{1mm}
  \dot{\xi}_3=-\dfrac{\partial H}{\partial x_3}=-\frac{1}{4}x_3(\theta_1^2+\theta_2^2-\theta_3^2)+\frac{1}{2}\xi_4\theta_1+\frac{1}{2}\xi_2\theta_2-\frac{1}{2}\xi_1\theta_3
  -\frac{1}{2}x_2\theta_1\theta_3+\frac{1}{2}x_4\theta_2\theta_3,  \\ \vspace{1mm}
  \dot{\xi}_4=-\dfrac{\partial H}{\partial x_4}=-\frac{1}{4}x_4(\theta_1^2-\theta_2^2+\theta_3^2)-\frac{1}{2}\xi_3\theta_1-\frac{1}{2}\xi_1\theta_2-\frac{1}{2}\xi_2\theta_3
  -\frac{1}{2}x_2\theta_1\theta_2+\frac{1}{2}x_3\theta_2\theta_3, \notag
\end{array}
\end{equation}
\begin{equation}
\begin{array}{l}
\hskip-7.5cm    \dot{\theta}_1=-\dfrac{\partial H}{\partial z_1}=0, \\ \vspace{2mm}
\hskip-7.5cm    \dot{\theta}_2=-\dfrac{\partial H}{\partial z_2}=0, \\ \vspace{1mm}
\hskip-7.5cm    \dot{\theta}_3=-\dfrac{\partial H}{\partial z_3}=0. \notag
\end{array}
\end{equation}

We observe that $\theta_1,\theta_2,\theta_3$ are constants.
Let us remind that the projection of a solution of the Hamiltonian
system onto $(x,z)$-space is called geodesic. In order to find it
we will reduce the Hamiltonian system  to the system containing
only $(x_1,x_2,x_3,x_4,z_1,z_2,z_3)$ coordinates. If we express
$\xi_1,\ldots,\xi_4$ from the first 4 equations and substitute them
in the equations of the Hamiltonian system, then we obtain
\begin{align}\notag
     &\dot{\xi}_1=\dfrac{1}{2}(\dot{x}_2\theta_1-\dot{x}_4\theta_2-\dot{x}_3\theta_3),\\ \vspace{1mm}\notag
     &\dot{\xi}_2=\dfrac{1}{2}(-\dot{x}_1\theta_1-\dot{x}_3\theta_2+\dot{x}_4\theta_3),\\ \vspace{1mm}\notag
     &\dot{\xi}_3=\dfrac{1}{2}(\dot{x}_4\theta_1+\dot{x}_2\theta_2+\dot{x}_1\theta_3),\\ \vspace{1mm}\notag
     &\dot{\xi}_4=\dfrac{1}{2}(-\dot{x}_3\theta_1+\dot{x}_1\theta_2-\dot{x}_2\theta_3).
\end{align}
Differentiating first 4 equations and substituting $\dot{\xi}_1,\ldots,\dot{\xi}_4$ there, we get
\begin{align*}
     &\ddot{x}_1=-\dot{x}_2\theta_1+\dot{x}_4\theta_2+\dot{x}_3\theta_3,\\
     &\ddot{x}_2=-\dot{x}_1\theta_1-\dot{x}_3\theta_2+\dot{x}_4\theta_3,\\
     &\ddot{x}_3=\dot{x}_4\theta_1+\dot{x}_2\theta_2+\dot{x}_1\theta_3,\\
     &\ddot{x}_4=-\dot{x}_3\theta_1+\dot{x}_1\theta_2-\dot{x}_2\theta_3
\end{align*}
or
\begin{equation}\label{sys:1}
    \left(%
   \begin{array}{l}
       \ddot{x}_1\\ \ddot{x}_2\\ \ddot{x}_3\\ \ddot{x}_4
   \end{array}%
    \right)=\left(%
   \begin{array}{cccc}
       0& -\theta_1 & \theta_3 & \theta_2 \\
      -\theta_1 & 0 &-\theta_2&\theta_3\\
      \theta_3&\theta_2&0&\theta_1\\
      \theta_2&-\theta_3&-\theta_1&0
   \end{array}%
    \right)\left(%
   \begin{array}{l}
       \dot{x}_1\\ \dot{x}_2\\ \dot{x}_3\\ \dot{x}_4
   \end{array}%
    \right).
\end{equation}
We are looking for the solution $x_1=x_1(t), \ldots,x_4=x_4(t)$, $t\in[-\infty,+\infty]$, satisfying $x_1(0)=0,\ldots,x_4(0)=0$ and $\dot{x}_1(0)=\dot{x}_1^0,\ldots,\dot{x}_4(0)=\dot{x}_4^0$.
The eigenvalues of the matrix
\begin{gather}
    A:=\left(%
   \begin{array}{cccc}
       0& -\theta_1 & \theta_3 & \theta_2 \\
      -\theta_1 & 0 &-\theta_2&\theta_3\\
      \theta_3&\theta_2&0&\theta_1\\
      \theta_2&-\theta_3&-\theta_1&0
   \end{array}%
    \right)\notag
\end{gather}
are  $\lambda_1=a$, $\lambda_2=-a$, $\lambda_3=ia$, and $\lambda_4=-ia$, where $a=\sqrt{\theta_1^2+\theta_2^2+\theta_3^2}$. The associated eigenvectors are
\begin{align*}
    &v_1=(a,-\theta_1,\theta_3,\theta_2),\\
    &v_2=(a,\theta_1,-\theta_3,-\theta_2),\\
    &v_3=(0,\theta_1\theta_3+ia\theta_2,\theta_1^2+\theta_2^2,ia\theta_1-\theta_2\theta_3),\\
    &v_4=(0,\theta_1\theta_3-ia\theta_2,\theta_1^2+\theta_2^2,-ia\theta_1-\theta_2\theta_3).
\end{align*}
The solution of the system (\ref{sys:1}) is of the form
\begin{align*}
    &\dot{x}_1=a(c_1+c_2)\cosh(at)+a(c_1-c_2)\sinh(at),\notag\\
    &\dot{x}_2=\theta_1(c_2-c_1)\cosh(at)-\theta_1(c_1+c_2)\sinh(at)\\
    &\hspace{7mm}+2(c_3\theta_1\theta_3-c_4a\theta_2)\cos(at)-2(c_4\theta_1\theta_3+c_3a\theta_2)\sin(at),\notag\\
    &\dot{x}_3= \theta_3(c_1-c_2)\cosh(at)+\theta_3(c_1+c_2)\sinh(at)\\
    &\hspace{7mm}+2c_3(\theta_1^2+\theta_2^2)\cos(at)-2c_4(\theta_1^2+\theta_2^2)\sin(at),\notag\\
    &\dot{x}_4= \theta_2(c_1-c_2)\cosh(at)+\theta_2(c_1+c_2)\sinh(at)\notag\\
    &\hspace{7mm}-2(c_4a\theta_1+c_3\theta_2\theta_3)\cos(at)-2(c_3a\theta_1-c_4\theta_2\theta_3)\sin(at),\notag
\end{align*}
where
\begin{align*}
    &c_1=\dfrac{a\dot{x}^0_1-\theta_1\dot{x}^0_2+\theta_3\dot{x}^0_3+\theta_2\dot{x}^0_4}{2a^2},\notag\\
    &c_2=\dfrac{a\dot{x}^0_1+\theta_1\dot{x}^0_2-\theta_3\dot{x}^0_3-\theta_2\dot{x}^0_4}{2a^2},\notag\\
    &c_3=\dfrac{\theta_1\theta_3\dot{x}^0_2+(\theta_1^2+\theta_2^2)\dot{x}^0_3-\theta_2\theta_3\dot{x}^0_4}{2a^2(\theta_1^2+\theta_2^2)},\notag\\
    &c_4=-\dfrac{\theta_2\dot{x}^0_2+\theta_1\dot{x}^0_4}{2a(\theta_1^2+\theta_2^2)}.\notag\\
\end{align*}
Therefore,
\begin{align}\label{eq:x}
    &x_1=A_1\sinh(at)+B_1\cosh(at)+E_1,\nonumber \\
    &x_2=A_2\sinh(at)+B_2\cosh(at)+C_2\sin(at)+D_2\cos(at)+E_2,\\
    &x_3=A_3\sinh(at)+B_3\cosh(at)+C_3\sin(at)+D_3\cos(at)+E_3,\nonumber \\
    &x_4=A_4\sinh(at)+B_4\cosh(at)+C_4\sin(at)+D_4\cos(at)+E_4,\nonumber
\end{align}
where the coefficients $A_i$, $B_i$, $C_j$, $D_j$, $E_i$ for
$i=1,2,3,4$, $j=1,2,3$ are given in terms of $c_i$, $\theta_j$ and are
presented in Appendix.

We find from the horizontality condition \eqref{eq:horizontality}  
\begin{gather}\label{eq:z}
    z_1=a\alpha^1_0 t+\beta^1_1 \sinh(at) \sin(at) +\beta^1_2 \cosh(at) \cos(at)\nonumber \\
\hspace{16mm}+\beta^1_3 \sinh(at) \cos(at)+\beta^1_4 \cosh(at)\sin(at)\nonumber\\
\hspace{46mm}    +\alpha^1_5\sinh(at)+\alpha^1_6\cosh(at)+\alpha^1_7\sin(at)+\alpha^1_8\cos(at)-\beta_2^1,
\nonumber \\
    z_2=a\alpha^2_0 t+\beta^2_1 \sinh(at) \sin(at) +\beta^2_2 \cosh(at) \cos(at) \\
\hspace{16mm}+\beta^2_3 \sinh(at) \cos(at)+\beta^2_4 \cosh(at)\sin(at)\nonumber \\
\hspace{46mm}     +\alpha^2_5\sinh(at)+\alpha^2_6\cosh(at)+\alpha^2_7\sin(at)+\alpha^2_8\cos(at)-\beta_2^2,
\nonumber \\
    z_3=a\alpha^3_0 t+\beta^3_1 \sinh(at) \sin(at) +\beta^3_2 \cosh(at) \cos(at)\nonumber \\
\hspace{16mm}+\beta^3_3 \sinh(at) \cos(at)+\beta^3_4 \cosh(at)\sin(at)\nonumber\\
\hspace{46mm}
+\alpha^3_5\sinh(at)+\alpha^3_6\cosh(at)+\alpha^3_7\sin(at)+\alpha^3_8\cos(at)-\beta_2^3.
\nonumber
\end{gather}
The coefficients $\alpha_{n}^{k}$ and $\beta_{m}^{k}$ can also be found
in Appendix. Then after long and laborious calculations we get
\begin{equation}\label{eq:x2}
-x_1^2+x_2^2+x_3^2+x_4^2=-16c_1c_2\sinh^2\bigl(\frac{at}{2}\bigr)+16k\sin^2\bigl(\frac{at}{2}\bigr).
\end{equation}
and
\begin{eqnarray}\label{eq:z2}
z_1^2+z_2^2+z_3^2=4(at(-c_1c_2+k)+2c_1c_2\sinh(at)-2k\sin(at))^2\notag\\-4k(c_1^2e^{at}+c_2^2e^{-at})(4\sin(at)\sinh(at)+5\cos(at)-5\cosh(at))
\notag\\
   +8c_1c_2k(5\sin(at)\sinh(at)+4\cos(at)-4\cosh(at)).
\end{eqnarray}
Notice that
$-c_1c_2+k=\frac{1}{4a^2}\|v_0\|^2$. The auxiliary calculations
can be found in Appendix.

We believe that the quantity of geodesics in case of quaternion
group equipped with the Lorentzian metric depends on the ratio
$\frac{z_1^2+z_2^2+z_3^2}{(x_1^2-x_2^2-x_3^2-x_4^2)^2}$ like in
the case of the Heisenberg group. The study of this question we
postpone for a forthcoming paper.

\section{Appendix}

The coefficients for the solutions~\eqref{eq:x}.
$$\begin{array}{lllll}
    & A_1=c_1+c_2,\qquad & A_2=\dfrac{\theta_1(c_2-c_1)}{a},\qquad & A_3=\dfrac{\theta_3(c_1-c_2)}{a},\qquad
    & A_4=\dfrac{\theta_2(c_1-c_2)}{a},\\
    & B_1=c_1-c_2, & B_2=-\dfrac{\theta_1(c_1+c_2)}{a}, & B_3=\dfrac{\theta_3(c_1+c_2)}{a},
    & B_4=\dfrac{\theta_2(c_1+c_2)}{a},
\end{array}$$
$$
\begin{array}{llll}
    & C_2=\dfrac{2(c_3\theta_1\theta_3-c_4a\theta_2)}{a},\quad &
    C_3=\dfrac{2c_3(\theta_1^2+\theta_2^2)}{a},\quad
    & C_4=-\dfrac{2(c_4a\theta_1+c_3\theta_2\theta_3)}{a},\\
    & D_2=\dfrac{2(c_4\theta_1\theta_3+c_3a\theta_2)}{a}, & D_3=\dfrac{2c_4(\theta_1^2+\theta_2^2)}{a},
    & D_4=\dfrac{2(c_3a\theta_1-c_4\theta_2\theta_3)}{a},
\end{array}
$$
$$
\begin{array}{lll}
    & E_1=c_2-c_1,  & E_2=\dfrac{\theta_1(c_1+c_2)}{a}-\dfrac{2(c_4\theta_1\theta_3+c_3a\theta_2)}{a},
    \\
    & E_3=-\dfrac{\theta_3(c_1+c_2)}{a}-\dfrac{2c_4(\theta_1^2+\theta_2^2)}{a}, \quad
    & E_4=-\dfrac{\theta_2(c_1+c_2)}{a}-\dfrac{2(c_3a\theta_1-c_4\theta_2\theta_3)}{a}.
\end{array}
$$

The coefficients $\alpha_{n}^{k}$ and $\beta_{m}^{k}$
in~\eqref{eq:z} are related in the following way.
$\beta^i_1=\dfrac{1}{4}(\alpha^i_1+\alpha^i_4)$,
$\beta^i_2=\dfrac{1}{4}(-\alpha^i_1+\alpha^i_4)$,
$\beta^i_3=\dfrac{1}{4}(-\alpha^i_2+\alpha^i_3)$,
$\beta^i_4=\dfrac{1}{4}(\alpha^i_2+\alpha^i_3)$, $i=1,\,2,\,3$.
The list of $\alpha_{n}^{k}$ is presented.
$$\begin{array}{lll}
\alpha^1_0 & = & \frac{2\theta_1}{a}(-c_1c_2+k), \\
\alpha^2_0 & = & \frac{2\theta_2}{a}(-c_1c_2+k), \\
\alpha^3_0 & = & \frac{2\theta_3}{a}(-c_1c_2+k),\\
\alpha^1_1 & = &
\frac{4}{a}(c_1+c_2)(c_3\theta_1\theta_3-ac_4\theta_2),\\
\alpha^2_1 & = &
\frac{4}{a}(c_1+c_2)(c_3\theta_2\theta_3+ac_4\theta_1),\\
\alpha^3_1 & = &
-\frac{4c_3}{a}(c_1+c_2)(\theta_1^2+\theta_2^2),\\
\end{array}$$
$$\begin{array}{lll}
\alpha^1_2 & = & \frac{4}{a}(c_1-c_2)(c_3\theta_1\theta_3-ac_4\theta_2),\\
\alpha^2_2 & = & \frac{4}{a}(c_1-c_2)(c_3\theta_2\theta_3+ac_4\theta_1),\\
\alpha^3_2 & = & -\frac{4c_3}{a}(c_1-c_2)(\theta_1^2+\theta_2^2),\\
\alpha^1_3 & = & \frac{4}{a}(c_1+c_2)(c_4\theta_1\theta_3+ac_3\theta_2),\\
\alpha^2_3 & = & \frac{4}{a}(c_1+c_2)(c_4\theta_2\theta_3-ac_3\theta_1),\\
\alpha^3_3 & = & -\frac{4c_4}{a}(c_1+c_2)(\theta_1^2+\theta_2^2),\\
\alpha^1_4 & = & \frac{4}{a}(c_1-c_2)(c_4\theta_1\theta_3+ac_3\theta_2),\\
\alpha^2_4 & = &\frac{4}{a}(c_1-c_2)(c_4\theta_2\theta_3-ac_3\theta_1),\\
\alpha^3_4 & = & -\frac{4c_4}{a}(c_1-c_2)(\theta_1^2+\theta_2^2),
\end{array}$$
$$
\begin{array}{ll}
   &\alpha^1_5= -2\theta_2(c_3(c_1+c_2)-c_4(c_1-c_2))+\frac{4c_1c_2\theta_1}{a}-\frac{2\theta_1\theta_3}{a}(c_3(c_1-c_2)+c_4(c_1+c_2)),\\
   &\alpha^2_5= 2\theta_1(c_3(c_1+c_2)-c_4(c_1-c_2))+\frac{4c_1c_2\theta_2}{a}-\frac{2\theta_2\theta_3}{a}(c_3(c_1-c_2)+c_4(c_1+c_2)),\\
   &\alpha^3_5= \frac{2}{a}(\theta_1^2+\theta_2^2)(c_3(c_1-c_2)+c_4(c_1+c_2))+\frac{4c_1c_2\theta_3}{a},\\
   &\alpha^1_6= 2\theta_2(-c_3(c_1-c_2)+c_4(c_1+c_2))-\frac{2\theta_1\theta_3}{a}(c_3(c_1+c_2)+c_4(c_1-c_2)),\\
   &\alpha^2_6= 2\theta_1(c_3(c_1-c_2)-c_4(c_1+c_2))-\frac{2\theta_2\theta_3}{a}(c_3(c_1+c_2)+c_4(c_1-c_2)),\\
   &\alpha^3_6= \frac{2}{a}(\theta_1^2+\theta_2^2)(c_3(c_1+c_2)+c_4(c_1-c_2)),\\
   &\alpha^1_7= -2\theta_2(c_3(c_1+c_2)+c_4(c_1-c_2))- \frac{2\theta_1\theta_3}{a}(-c_3(c_1-c_2)+c_4(c_1+c_2))-\frac{4}{a}\theta_1k,\\
   &\alpha^2_7= 2\theta_1(c_3(c_1+c_2)+c_4(c_1-c_2))- \frac{2\theta_2\theta_3}{a}(-c_3(c_1-c_2)+c_4(c_1+c_2))-\frac{4}{a}\theta_2k,\\
   &\alpha^3_7= -\frac{2}{a}(\theta_1^2+\theta_2^2)(c_3(c_1-c_2)-c_4(c_1+c_2))-\frac{4}{a}\theta_3k,\\
   &\alpha^1_8= -\alpha_6^1,\quad
   \alpha^2_8=-\alpha^2_6,\quad
   \alpha^3_8= -\alpha^3_6,
\end{array}$$ where $k=(c_3^2+c_4^2)(\theta_1^2+\theta_2^2)$.

In the following lines we give some useful calculations to
deduce~\eqref{eq:x2}.
$$
\begin{array}{ccc}
 -A_1^2+A_2^2+A_3^2+A_4^2=-4c_1c_2,\quad -B_1^2+B_2^2+B_3^2+B_4^2=4c_1c_2,\notag \\
C_2^2+C_3^2+C_4^2=D_2^2+D_3^2+D_4^2=4k,\quad -E_1^2+E_2^2+E_3^2+E_4^2=4(c_1c_2+k),
\end{array}
$$
$$
-A_1B_1+A_2B_2+A_3B_3+A_4B_4=
 A_2C_2+A_3C_3+A_4C_4=A_2D_2+A_3D_3+A_4D_4=0,
$$
$$
-A_1E_1+A_2E_2+A_3E_3+A_4E_4=0,
$$
$$
B_2C_2+B_3C_3+B_4C_4= B_2D_2+B_3D_3+B_4D_4=0,
$$
$$
-B_1E_1+B_2E_2+B_3E_3+B_4E_4=-4c_1c_2,
$$
$$
C_2D_2+C_3D_3+C_4D_4= C_2E_3+C_3E_3+C_4E_4=0,\quad
D_2E_2+D_3E_3+D_4E_4=-4k.
$$
Now we present the values of different combinations that we used
in calculation of $z_1^2+z_2^2+z_3^2$.

$$
   (\alpha_0^1)^2+(\alpha_0^2)^2+(\alpha_0^3)^2=4(-c_1c_2+k)^2
   =\frac{1}{4a^4}(-(\dot{x}_1^0)^2+(\dot{x}_2^0)^2+(\dot{x}_3^0)^2+(\dot{x}_4^0)^2)^2,
$$
\begin{gather*}
   (\beta_1^1)^2+(\beta_1^2)^2+(\beta_1^3)^2=(\beta_2^1)^2+(\beta_2^2)^2+(\beta_2^3)^2=(\beta_3^1)^2+(\beta_3^2)^2+(\beta_3^3)^2\notag\\
   =(\beta_4^1)^2+(\beta_4^2)^2+(\beta_4^3)^2=2(c_1^2+c_2^2)k,
\end{gather*}
$$
\begin{array}{ll}
&  (\alpha_5^1)^2+(\alpha_5^2)^2+(\alpha_5^3)^2  =  8(c_1^2+c_2^2)k+16c_1^2c_2^2,\notag\\
&  (\alpha_6^1)^2+(\alpha_6^2)^2+(\alpha_6^3)^2  =  8(c_1^2+c_2^2)k,\notag\\
&   (\alpha_7^1)^2+(\alpha_7^2)^2+(\alpha_7^3)^2  =  8(c_1^2+c_2^2)k+16k^2,\notag\\
&   (\alpha_8^1)^2+(\alpha_8^2)^2+(\alpha_8^3)^2  =
8(c_1^2+c_2^2)k,
\end{array}
$$
$$
\alpha_0^1\beta_1^1+\alpha_0^2\beta_1^2+\alpha_0^3\beta_1^3
=\alpha_0^1\beta_2^1+\alpha_0^2\beta_2^2+\alpha_0^3\beta_2^3=
   \alpha_0^1\beta_3^1+\alpha_0^2\beta_3^2+\alpha_0^3\beta_3^3
   =\alpha_0^1\beta_4^1+\alpha_0^2\beta_4^2+\alpha_0^3\beta_4^3=0,
$$
$$
\begin{array}{lll}
&
\alpha_0^1\alpha_5^1+\alpha_0^2\alpha_5^2+\alpha_0^3\alpha_5^3=-8c_1^2c_2^2+8c_1c_2k,\quad
&
   \alpha_0^1\alpha_6^1+\alpha_0^2\alpha_6^2+\alpha_0^3\alpha_6^3=0,\notag\\
&
\alpha_0^1\alpha_7^1+\alpha_0^2\alpha_7^2+\alpha_0^3\alpha_7^3=8c_1c_2k-8k^2,
&
\alpha_0^1\alpha_8^1+\alpha_0^2\alpha_8^2+\alpha_0^3\alpha_8^3=0,
\end{array}
$$
\begin{gather*}
\beta_1^1\beta_2^1+\beta_3^1\beta_4^1+\beta_1^2\beta_2^2
+\beta_3^2\beta_4^2+\beta_1^3\beta_2^3+\beta_3^3\beta_4^3=0,\notag\\
   \beta_1^1\beta_3^1+\beta_1^2\beta_3^2+\beta_1^3\beta_3^3=0,\notag\\
   \beta_1^1\beta_4^1+\beta_1^2\beta_4^2+\beta_1^3\beta_4^3=2(c_1^2-c_2^2)k,\notag\\
   \alpha_5^1\beta_1^1+\alpha_5^2\beta_1^2+\alpha_5^3\beta_1^3=-4(c_1^2-c_2^2)k,\notag\\
   \alpha_6^1\beta_1^1+\alpha_5^1\beta_4^1+\alpha_6^2\beta_1^2+\alpha_5^2\beta_4^2+\alpha_6^3\beta_1^3+\alpha_5^3\beta_4^3=-8(c_1^2+c_2^2)k,
\end{gather*}
$$
\begin{array}{lll}
  & \alpha_7^1\beta_1^1+\alpha_7^2\beta_1^2+\alpha_7^3\beta_1^3=0,\quad
  & \alpha_7^1\beta_3^1+\alpha_8^1\beta_1^1
   +\alpha_7^2\beta_3^2+\alpha_8^2\beta_1^2+\alpha_7^3\beta_3^3+\alpha_8^3\beta_1^3=0,\notag\\
  &
  \beta_2^1\beta_3^1+\beta_2^2\beta_3^2+\beta_2^3\beta_3^3=2(c_1^2-c_2^2)k,\qquad
  & \beta_2^1\beta_4^1+\beta_2^2\beta_4^2+\beta_2^3\beta_4^3=0,
\end{array}
$$
\begin{gather*}
   \alpha_5^1\beta_2^1+\alpha_6^1\beta_3^1+\alpha_5^2\beta_2^2+\alpha_6^2\beta_3^2
   +\alpha_5^3\beta_2^3+\alpha_6^3\beta_3^3=0,\notag\\
   \alpha_6^1\beta_2^1+\alpha_6^2\beta_2^2+\alpha_6^3\beta_2^3=8c_1c_2k,\notag\\
   \alpha_7^1\beta_2^1+\alpha_8^1\beta_4^1+\alpha_7^2\beta_2^2+\alpha_8^2\beta_4^2
   +\alpha_7^3\beta_2^3+\alpha_8^3\beta_4^3=0,
\end{gather*}
$$
\begin{array}{lll}
&
\alpha_8^1\beta_2^1+\alpha_8^2\beta_2^2+\alpha_8^3\beta_2^3=-8c_1c_2k,\quad
&
   \alpha_5^1\beta_3^1+\alpha_5^2\beta_3^2+\alpha_5^3\beta_3^3=-8c_1c_2k,\notag\\
&
\alpha_8^1\beta_3^1+\alpha_8^2\beta_3^2+\alpha_8^3\beta_3^3=0,\quad
&
\alpha_6^1\beta_4^1+\alpha_6^2\beta_4^2+\alpha_6^3\beta_4^3=-4(c_1^2-c_2^2)k,
\end{array}
$$
$$
   \alpha_7^1\beta_4^1+\alpha_7^2\beta_4^2+\alpha_7^3\beta_4^3=-8c_1c_2k,
$$
$$
\begin{array}{lll}
&
\alpha_5^1\alpha_6^1+\alpha_5^2\alpha_6^2+\alpha_5^3\alpha_6^3=8(c_1^2-c_2^2)k,\quad
&   \alpha_5^1\alpha_7^1+\alpha_5^2\alpha_7^2+\alpha_5^3\alpha_7^3=0,\notag\\
&
\alpha_5^1\alpha_8^1+\alpha_5^2\alpha_8^2+\alpha_5^3\alpha_8^3=-8(c_1^2-c_2^2)k,\quad
&   \alpha_6^1\alpha_7^1+\alpha_6^2\alpha_7^2+\alpha_6^3\alpha_7^3=0,\notag\\
&
\alpha_6^1\alpha_8^1+\alpha_6^2\alpha_8^2+\alpha_6^3\alpha_8^3=-8(c_1^2+c_2^2)k,\quad
&   \alpha_7^1\alpha_8^1+\alpha_7^2\alpha_8^2+\alpha_7^3\alpha_8^3=0,
\end{array}
$$
$$
   \beta_1^1\beta_2^1+\beta_1^2\beta_2^2+\beta_1^3\beta_2^3=-4c_1c_2k,
$$
$$
   \alpha_5^1\beta_2^1+\alpha_5^2\beta_2^2+\alpha_5^3\beta_2^3=0,\quad
   \alpha_7^1\beta_2^1+\alpha_7^2\beta_2^2+\alpha_7^3\beta_2^3=-4(c_1^2-c_2^2)k,\quad
   \alpha_8^1\beta_2^1+\alpha_8^2\beta_2^2+\alpha_8^3\beta_2^3=-8c_1c_2k.
$$

\end{document}